\documentclass{article}

\usepackage{arxiv}

\usepackage[T1]{fontenc}
\usepackage[utf8]{inputenc}
\usepackage{times}

\usepackage{amsmath}
\usepackage{amssymb}
\usepackage{amsthm}
\numberwithin{equation}{section}

\usepackage[ngerman,english]{babel}

\usepackage{graphicx}
\usepackage{booktabs}
\usepackage{caption}
\usepackage{subcaption}
\usepackage{makecell}

\usepackage{cite}
\usepackage{hyperref}
\usepackage{doi}

\usepackage{ulem}

\usepackage{ifthen}
\newcommand{\mycomment}[1]{%
  \ifthenelse{\isodd{\value{page}}}{%
    \reversemarginpar\marginpar{\tiny{#1}}%
  }{%
    \reversemarginpar\marginpar{\tiny{#1}}%
  }}

\newcommand{\takeout}[1]{}

\title{Four-Level Overlapping Schwarz as Multigrid Coarse Solver for Incompressible Non-Newtonian Flow in Complex Geometries}

\author{
  Stephan Köhler\thanks{%
    \texttt{stephan.koehler@math.tu-freiberg.de},
    phone +49\,3731\,39\,3188} \\
  INMO, Technische Universität Bergakademie Freiberg \\
  09596 Freiberg, Germany
  \And
  Oliver Rheinbach\thanks{%
    \texttt{oliver.rheinbach@math.tu-freiberg.de},
    phone +49\,3731\,39\,3279} \\
  INMO / URZ / ZeHS, Technische Universität Bergakademie Freiberg \\
  09596 Freiberg, Germany
}

\newenvironment{acknowledgement}{%
  \section*{Acknowledgements}%
}{}

\date{}   

\newenvironment{vchtable}[1][tbp]{%
  \begin{table}[#1]%
  }{%
  \end{table}%
}

\begin{document}
\maketitle

\begin{abstract}
  For complex geometries, the coarse problem of geometric multigrid can be too large to be solved by a direct solver.
  Here, we report on the use of domain decomposition applied to the multigrid coarse problem.
  Additive overlapping Schwarz methods are domain decomposition methods for the iterative solution of partial differential equations whose numerical and parallel scalability can be improved by the addition of coarse levels. A successful coarse space for such methods, inspired by iterative substructuring, is the generalized Dryja--Smith--Widlund (GDSW) space. A monolithic two-level overlapping Schwarz preconditioner based on a GDSW coarse space has been introduced for the solution of saddle-point problems arising from incompressible fluid problems, and has subsequently been extended to a three-level method.
  In the present work, for the first time, we consider a monolithic four-level overlapping Schwarz preconditioner, obtained by applying the two-level monolithic GDSW construction recursively three times, so that the second- and third-level coarse problems are themselves treated by overlapping Schwarz and only the smallest fourth-level coarse problem is solved by a sparse direct method.
  The preconditioner is implemented in the FROSch (Fast and Robust Overlapping Schwarz) library, part of the Trilinos package ShyLU, and is coupled to the FEATFLOW library through a dedicated scalable interface that makes FROSch available both as a preconditioner for the full saddle-point system and as a coarse solver inside the FEATFLOW geometric multigrid method. Numerical results are presented for a three-dimensional incompressible stationary Stokes problem with a Carreau-type non-Newtonian viscosity model posed on the complex geometry of an extrusion die, on up to $4\,000$\,MPI ranks, comparing the four-level preconditioner with its two-level and three-level counterparts in both roles. This work is part of the StroemungsRaum project, funded by the German Bundesministerium für Forschung, Technologie und Raumfahrt (BMFTR, formerly BMBF) as part of the SCALEXA program on new methods and technologies for exascale computing.
\end{abstract}

\maketitle

\section{Introduction}
The discretization of fluid problems of Stokes or Navier--Stokes type leads to linear systems of saddle-point structure. As the geometric complexity of the domain increases and finer meshes are required to resolve boundary layers or other small-scale flow features, these systems quickly grow to a size at which sparse direct solvers are no longer a viable option, both in terms of arithmetic cost and in terms of memory consumption. This makes Krylov subspace methods, equipped with strong preconditioners, the method of choice. Among such preconditioners, multilevel approaches from multigrid and domain decomposition are particularly attractive on modern HPC architectures, since they distribute both the solution effort and the associated data across many compute nodes while still providing a good convergence rate.

We consider monolithic overlapping Schwarz preconditioners of the generalized Dryja--Smith--Widlund (GDSW) type~\cite{Dohrmann:2008:FEM,Dohrmann:2008:DDL} for incompressible fluid problems. In the monolithic formulation, the local overlapping problems and the coarse problem retain the saddle-point structure of the original system; no static condensation of the pressure is performed. This is the framework introduced in~\cite{Heinlein:2019:Monolithic}; see also~\cite{Hochmuth:2020:Parallel}. While the two-level variant is already numerically scalable for moderate subdomain counts, its parallel scalability is fundamentally limited by the size of the GDSW coarse problem, which grows with the number of subdomains. Once this coarse problem can no longer be factorized by a sparse direct solver, the two-level method ceases to be applicable.

A standard remedy is to extend the construction by additional coarse levels. Applying the two-level idea recursively to the GDSW coarse problem yields the three-level overlapping Schwarz preconditioner of~\cite{Heinlein:2018:Three,Roever:2019:ThreeLevel}; see also~\cite{Roever:2022:Multilevel:PHD,Heinlein:2023:Multilevel}. The three-level method substantially enlarges the range of tractable subdomain counts and, once the coarse problem is sufficiently large, is also faster than the two-level method, because the cost of the direct factorization of the coarse problem dominates in the two-level case.
However, the same scalability barrier reappears, one level deeper, as soon as the three-level coarse problem itself grows large: its factorization then becomes the new bottleneck.

For the first time, we present results for a monolithic four-level overlapping Schwarz preconditioner, in which the two-level monolithic GDSW construction is applied recursively three times, so that the second- and the third-level coarse problems are themselves treated by overlapping Schwarz, and only the smallest, fourth-level coarse problem is solved by a sparse direct solver.

We show that our monolithic three-level overlapping Schwarz preconditioner, implemented in the ShyLU/FROSch package~\cite{Shylu:2012,Heinlein:2018:FPI}, can be extended in practice to a four-level variant without any extra implementational overhead. The monolithic two-level overlapping Schwarz preconditioner with a GDSW type coarse space was introduced in~\cite{Heinlein:2019:Monolithic} and extended to a multi-level implementation in~\cite{Koehler:2026:Monolithic,Koehler:2026:Combining} based on the multi-level extension introduced in~\cite{Roever:2019:ThreeLevel}.

To assess the resulting method, we present numerical results up to 4\,000\,MPI ranks for a three-dimensional flow problem on an extrusion die geometry with a highly viscous non-Newtonian Fluid. The experiments compare the monolithic four-level preconditioner with the two-level and three-level variants.

This work is embedded in the StroemungsRaum project, funded by the BMFTR (formerly BMBF) under the SCALEXA program, the German initiative for software development in the exascale era. StroemungsRaum extends the CFD software FEATFLOW~\cite{featflow-website} by highly scalable domain decomposition methods and multigrid solvers tailored for heterogeneous CPU--GPU systems, with the goal of efficient CFD simulations on emerging exascale architectures. FEATFLOW is a modern C++ finite element software developed at TU Dortmund as a modular successor to earlier members of the FEAT family. It emphasizes highly scalable geometric multigrid solvers and is used within the project together with the industrial partner IANUS Simulation, who applies it for simulation-as-a-service solutions.
Within StroemungsRaum, the overlapping Schwarz preconditioner is used when the problem on the coarsest level of geometric multigrid is, due to complex geometries, too large to be solved directly.

The remainder of the paper is organized as follows. In Section~\ref{pamm-26-koehler-sec:stokes-problem}, the model Stokes problem is introduced together with a brief summary of the monolithic overlapping Schwarz preconditioner, followed by the extension to the monolithic multi-level case and, in particular, to the four-level construction. Numerical results are reported in Section~\ref{pamm-26-koehler-sec:model-problem}.

\section{Stokes Problem}\label{pamm-26-koehler-sec:stokes-problem}
Let $\Omega\subset\mathbf{R}^{d}$, $d\in\{2,3\}$, be a bounded Lipschitz domain whose boundary is partitioned into a Dirichlet part $\partial\Omega_{D}$ and a Neumann part $\partial\Omega_{N}$, with $n$ denoting the outward unit normal. We consider the stationary Stokes problem with viscosity $\mu>0$, body force $f$ and Dirichlet datum $g$. The weak formulation is given by: Find a velocity $u\in V_{g}$ and a pressure $p\in Q$ such that
\begin{equation*}
  \begin{aligned}
    \mu \int_{\Omega}\nabla u\colon\nabla v\,\mathrm{d}x &\,\, - \int_{\Omega}\mathrm{div}v\,p\,\mathrm{d}x &=&  \int_{\Omega}fv\mathrm{d}x & & \forall v\in V_{0}, \\
    -\int_{\Omega}\mathrm{div}u\,q\,\mathrm{d}x          &                                                 &=&\,  0                      & & \forall q\in Q_{0}, \\[1.5ex]
    u = g, &  \quad x\in\partial\Omega_{D}, & &&  \frac{\partial u}{\partial n} -pn = 0,  &\quad x\in\partial\Omega_{N}.
  \end{aligned}
\end{equation*}
Since all configurations considered in the numerical section combine an inflow and an outflow boundary, the pressure is uniquely determined; the non-unique case can be handled as in~\cite{Heinlein:2019:Monolithic}. We discretize with the FEATFLOW library~\cite{featflow-website} using $Q_{2}$ elements for the velocity and $P_{1}$-discontinuous elements for the pressure. The resulting algebraic system has the saddle-point form
\begin{equation}\label{eq:pamm-feat-frosch-saddle-point}
  \mathcal{A}\begin{bmatrix} u \\[0.5ex] p \end{bmatrix}
  = \mathcal{F},
  \qquad
  \mathcal{A} := \begin{bmatrix} A & B^{T} \\[0.5ex] B & 0 \end{bmatrix},
  \qquad
  \mathcal{F} := \begin{bmatrix} F \\[0.5ex] 0 \end{bmatrix},
\end{equation}
where $A$ corresponds to the discrete velocity Laplacian and $B^{T}$ couples velocity and pressure.

\subsection{Monolithic Overlapping Schwarz Preconditioner}

The construction of a monolithic overlapping Schwarz preconditioner for~\eqref{eq:pamm-feat-frosch-saddle-point} proceeds in two stages: an overlapping decomposition supplies the first level, and a GDSW coarse space supplies the second.

For the first level, the computational domain $\Omega$ is partitioned into $N$ nonoverlapping subdomains $\Omega_{i}$, $i=1,\ldots,N$, and each $\Omega_{i}$ is enlarged by $k$ layers of elements to obtain the overlapping subdomain $\Omega_{i}'$. The block restriction operator and the local saddle-point problem on $\Omega_{i}'$ are
\begin{equation*}
  \mathcal{R}_{i} := \begin{bmatrix} R_{i,u} & 0 \\[0.5ex] 0 & R_{i,p} \end{bmatrix},
  \qquad
  \mathcal{A}_{i} := \mathcal{R}_{i}\,\mathcal{A}\,\mathcal{R}_{i}^{T},
\end{equation*}
where $R_{i,u}$ and $R_{i,p}$ are the restrictions of the velocity and the pressure to $\Omega_{i}'$, respectively. Each $\mathcal{A}_{i}$ corresponds to a Stokes-type problem with homogeneous Dirichlet data on $\partial\Omega_{i}'\setminus\partial\Omega_{D}$ and is therefore invertible.

For the second level, the GDSW coarse space, we define the interface as $\Gamma := \bigcup_{i=1}^{N}\partial\Omega_{i}\setminus\partial\Omega_{D}$ and reorder the degrees of freedom (dofs) so that interface dofs (subscript $\Gamma$) are separated from the interior dofs of each subdomain (subscript $I$):
\begin{equation*}
  \mathcal{A} = \begin{bmatrix} \mathcal{A}_{II} & \mathcal{A}_{I\Gamma} \\[0.5ex] \mathcal{A}_{\Gamma I} & \mathcal{A}_{\Gamma\Gamma} \end{bmatrix}.
\end{equation*}
The interface is further decomposed into vertices, edges and, in 3D, faces. The columns of the coarse basis matrix $\Phi$ are obtained by prescribing nodal values on the interface entities according to the (rigid-body) null space of the differential operators and then extending these traces discretely as a saddle-point extension into the subdomain interiors:
\begin{equation*}
  \Phi := \begin{bmatrix} \Phi_{I} \\[0.5ex] \Phi_{\Gamma} \end{bmatrix}
        = \begin{bmatrix} -\mathcal{A}_{II}^{-1}\mathcal{A}_{I\Gamma}\Phi_{\Gamma} \\[0.5ex] \Phi_{\Gamma} \end{bmatrix}.
\end{equation*}
Each column of the interface part $\Phi_{\Gamma}$ corresponds to a single interface entity (vertex, edge, or face) and to a single null-space mode, and its entries are the values of this mode at the dofs of that entity.

Since the saddle-point problem~\eqref{eq:pamm-feat-frosch-saddle-point} carries separate null spaces for the velocity (translations and, optionally, rotations) and for the pressure (constants), the interface part of the basis is block-diagonal,
\begin{equation*}
  \Phi_{\Gamma} = \begin{bmatrix} \Phi_{\Gamma, u_{0}} & 0 \\[0.5ex] 0 & \Phi_{\Gamma, p_{0}} \end{bmatrix},
\end{equation*}
whereas the interior part of $\Phi$ is, in general, fully coupled across velocity and pressure:
\begin{equation*}
  \Phi
  =
  \begin{bmatrix}
    \Phi_{I,u_{0}, u} & \Phi_{I,u_{0}, p} \\[0.5ex]
    \Phi_{I,p_{0}, u} & \Phi_{I,p_{0}, p} \\[0.5ex]
    \Phi_{\Gamma, u_{0}} & 0 \\[0.5ex]
    0 & \Phi_{\Gamma, p_{0}}
  \end{bmatrix}.
\end{equation*}
For a discussion of the off-diagonal blocks $\Phi_{I,u_{0}, p}$ and $\Phi_{I,p_{0}, u}$ we refer to~\cite{Heinlein:2019:Monolithic,Hochmuth:2020:Parallel}. As in the original construction, only translations are used for the velocity null space.

With this, the monolithic two-level GDSW coarse operator and the resulting preconditioner read
\begin{equation*}
  \mathcal{A}_{0} := \Phi^{T}\mathcal{A}\,\Phi,
  \qquad
  \mathcal{B}_{\mathrm{GDSW}}
  = \Phi\,\mathcal{A}_{0}^{-1}\Phi^{T}
  + \sum_{i=1}^{N}\mathcal{R}_{i}^{T}\mathcal{A}_{i}^{-1}\mathcal{R}_{i}.
\end{equation*}
The coarse operator $\mathcal{A}_{0}$ is invertible due to the unique pressure induced by the outflow boundary condition. We remark that the GDSW preconditioner can be constructed in an almost entirely algebraic manner. Two exceptions are noted. First, the use of rotations in the velocity null space cannot be inferred from the matrix alone and is therefore not employed here. Second, the coupling among the pressure degrees of freedom cannot be extracted from $\mathcal{A}_{0}$ and likewise requires additional information. Since we use P1-discontinuous elements for the pressure, we know this coupling and can provide this information to FROSch.

\subsection{Monolithic Multi-Level Overlapping Schwarz Preconditioner}

As the number of subdomains $N$ grows, the coarse problem $\mathcal{A}_{0}$ grows accordingly, and its direct factorization becomes prohibitive in both time and memory. The multi-level remedy proposed in~\cite{Heinlein:2018:Three,Roever:2022:Multilevel:PHD,Heinlein:2023:Multilevel} is to apply the two-level construction recursively: the action of $\mathcal{A}_{0}^{-1}$ is itself approximated by a monolithic two-level overlapping Schwarz preconditioner built on a coarser decomposition. Repeating this step yields a hierarchy of progressively smaller coarse problems. In the present work, this recursion is carried out three times, which gives the monolithic four-level preconditioner.

To formalize this, we introduce a level index $\ell\in\{1,2,3,4\}$, where $\ell=1$ is the finest level and $\ell=4$ is the coarsest. The coarse operators are defined recursively by
\begin{equation*}
  \mathcal{A}_{0}^{(1)} := \mathcal{A},
  \qquad
  \mathcal{A}_{0}^{(\ell)} := \left.\Phi^{(\ell)}\right.^{T}\mathcal{A}_{0}^{(\ell-1)}\,\Phi^{(\ell)},
  \qquad \ell = 2, 3, 4,
\end{equation*}
where the columns of $\Phi^{(\ell)}$ are the monolithic GDSW coarse basis functions associated with an overlapping decomposition of the level-$(\ell-1)$ domain into $N^{(\ell-1)}$ subdomains. For $\ell\in\{1,2,3\}$, the level-$\ell$ block restrictions and the local saddle-point operators are denoted by $\mathcal{R}_{i}^{(\ell)}$ and $\mathcal{A}_{i}^{(\ell)}$, $i=1,\ldots,N^{(\ell)}$, respectively. Here, $N^{(1)}$ counts the (finest) subdomains, $N^{(2)}$ the subregions, and $N^{(3)}$ the subregions.

With this notation, the monolithic four-level overlapping Schwarz preconditioner is given by
\begin{equation}\label{eq:pamm-feat-frosch-gdsw-four-level}
  \begin{aligned}
    &\mathcal{B}_{\mathrm{GDSW},4\text{-}\mathrm{level}}\\[0.5ex]
    =\ &
    \underbrace{\Phi^{(2)}\!\left(
      \underbrace{\Phi^{(3)}\!\left(\underbrace{
        \underbrace{\Phi^{(4)}\left.\mathcal{A}_{0}^{(4)}\right.^{-1}\left.\Phi^{(4)}\right.^{T}}_{\text{4th level}}
        + \underbrace{\sum_{i=1}^{N^{(3)}}\left.\mathcal{R}_{i}^{(3)}\right.^{T}\left.\mathcal{A}_{i}^{(3)}\right.^{-1}\mathcal{R}_{i}^{(3)}}_{\text{3rd level}}
         }_{\text{approximation of } \left.\mathcal{A}_{0}^{(3)}\right.^{-1}}\right)\left.\Phi^{(3)}\right.^{T}
      + \underbrace{\sum_{i=1}^{N^{(2)}}\left.\mathcal{R}_{i}^{(2)}\right.^{T}\left.\mathcal{A}_{i}^{(2)}\right.^{-1}\mathcal{R}_{i}^{(2)}}_{\text{2nd level}}
    }_{\text{approximation of } \left.\mathcal{A}_{0}^{(2)}\right.^{-1}}
    \right)\left.\Phi^{(2)}\right.^{T}}_{\text{coarse levels}} \\
    & + \underbrace{\sum_{i=1}^{N^{(1)}}\left.\mathcal{R}_{i}^{(1)}\right.^{T}\left.\mathcal{A}_{i}^{(1)}\right.^{-1}\mathcal{R}_{i}^{(1)}}_{\text{1st level}}.
  \end{aligned}
\end{equation}
The first-level sum is the standard one-level overlapping Schwarz contribution on the original problem. The second- and third-level sums are the analogous contributions on the corresponding coarse decompositions, and only the fourth-level coarse problem $\mathcal{A}_{0}^{(4)}$ is solved by a direct method.

\section{Model problem}\label{pamm-26-koehler-sec:model-problem}

\begin{figure}[tbh]
  \centering
  \includegraphics[width=0.35\textwidth]{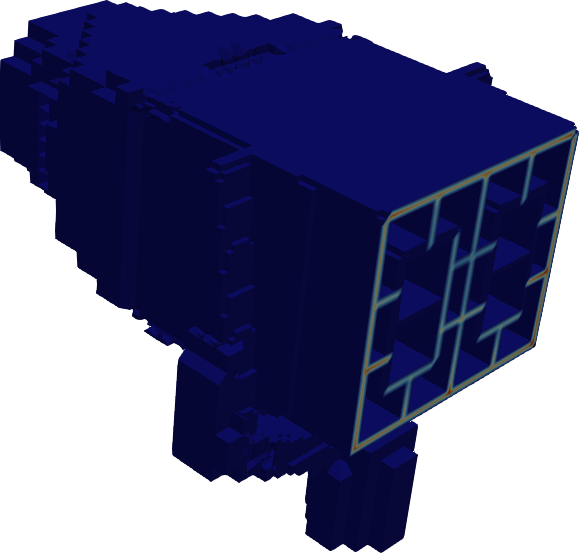} \hspace*{5ex}
  \includegraphics[width=0.40\textwidth]{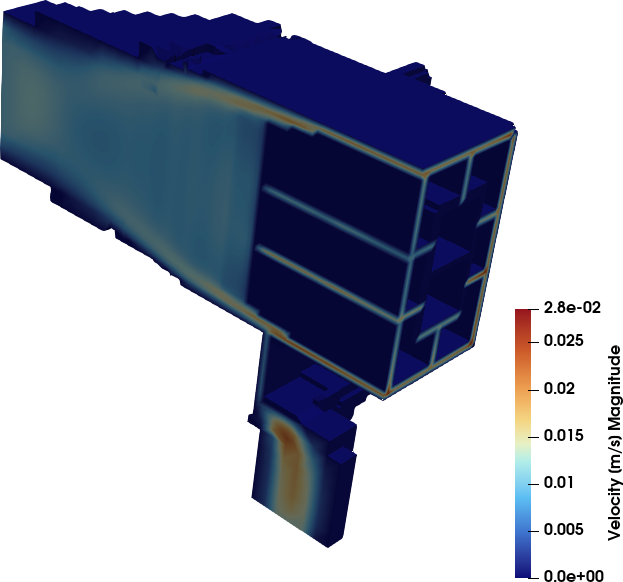}

  \caption{Visualization of the magnitude of the velocity of an incompressible Navier--Stokes simulation with a Carreau-type material model on an extrusion die. There are two inflow boundaries and one outflow boundary.}\label{pamm-26-koehler-fig:gendie}
\end{figure}

As a model problem, we consider the stationary flow of a highly viscous, incompressible non-Newtonian fluid governed by a Stokes-type system with a Carreau constitutive law~\cite{Carreau:1972:Rheological}. For an extrusion die, see Fig.~\ref{pamm-26-koehler-fig:gendie}, with boundary partitioned into an inflow part $\partial\Omega_{\mathrm{in}}$, an outflow part $\partial\Omega_{\mathrm{out}}$, and a no-slip part $\partial\Omega_{\mathrm{wall}} := \partial\Omega\setminus(\partial\Omega_{\mathrm{in}}\cup\partial\Omega_{\mathrm{out}})$, we seek a velocity $\mathbf{u}\in H^{1}(\Omega)^{d}$ and a pressure $p\in L^{2}(\Omega)$ satisfying
\begin{align*}
  -\nabla p + \nabla\cdot\boldsymbol{\sigma}(\mathbf{u},p) &= 0 && \text{in } \Omega, \\
  \nabla\cdot\mathbf{u} &= 0 && \text{in } \Omega, \\
  \mathbf{u} &= \mathbf{g} && \text{on } \partial\Omega_{\mathrm{in}}, \\
  \boldsymbol{\sigma}(\mathbf{u},p)\,\mathbf{n} &= 0 && \text{on } \partial\Omega_{\mathrm{out}}, \\
  \mathbf{u} &= 0 && \text{on } \partial\Omega_{\mathrm{wall}}.
\end{align*}
The constitutive law is given by
\begin{equation*}
  \boldsymbol{\sigma}(\mathbf{u},p) = -p\,\mathbf{I} + 2\,\nu(\mathbf{u})\,\mathbf{D}(\mathbf{u}),
  \qquad
  \mathbf{D}(\mathbf{u}) := \tfrac{1}{2}\!\left(\nabla\mathbf{u} + (\nabla\mathbf{u})^{T}\right),
\end{equation*}
with a shear-rate-dependent viscosity of Carreau type,
\begin{equation*}
  \nu(\mathbf{u}) = \nu_{0}\!\left(1 + (\lambda\,\dot{\gamma}(\mathbf{u}))^{2}\right)^{\!\frac{n-1}{2}},
  \qquad
  \dot{\gamma}(\mathbf{u}) := \sqrt{2\,\mathbf{D}(\mathbf{u})\colon\mathbf{D}(\mathbf{u})},
\end{equation*}
and a high zero-shear viscosity $\nu_{0}\gg 1$ that places the problem in a highly viscous regime. The nonlinearity of the system enters exclusively through the viscosity $\nu(\mathbf{u})$.

The system is discretized by $Q_{2}$ finite elements for the velocity and $P_{1}$-discontinuous elements for the pressure. We use one layer of overlap for all levels. The computational domain $\Omega$ is the extrusion die shown in Fig.~\ref{pamm-26-koehler-fig:gendie}; for the engineering background of extrusion dies we refer to~\cite{Hopmann:2016:Extrusion}. The discretized problem is solved with the FEATFLOW library using a geometric multigrid method~\cite{munster:2012:FEM-BEM,Turek:2003:Fictitious}.

We present numerical results for two different approaches to the linearized system arising in each step of the outer nonlinear iteration. In the first approach, the linearized system is solved directly by a (restarted) FGMRES method preconditioned by our parallel monolithic multilevel FROSch implementation~\cite{Koehler:2026:Monolithic}. In the second approach, the linearized system is solved by a geometric multigrid method in which the coarsest-level problem is treated by a (restarted) FGMRES iteration preconditioned by the same monolithic multilevel FROSch implementation.

In both approaches, the outer nonlinear iteration is an alternating Picard--Newton scheme, which arises as the limit of an implicit-explicit method~\cite{Barrenechea:2024:Implicit}. In the multigrid approach, the coarsest-level problem is solved only to low accuracy, with a relative residual reduction of $10^{-1}$ as stopping criterion for the FROSch-preconditioned FGMRES. All computations were carried out on the HPC cluster of the Technische Universität Bergakademie Freiberg.

\begin{figure}[tbh]
  \centering
  \includegraphics[width=0.40\textwidth]{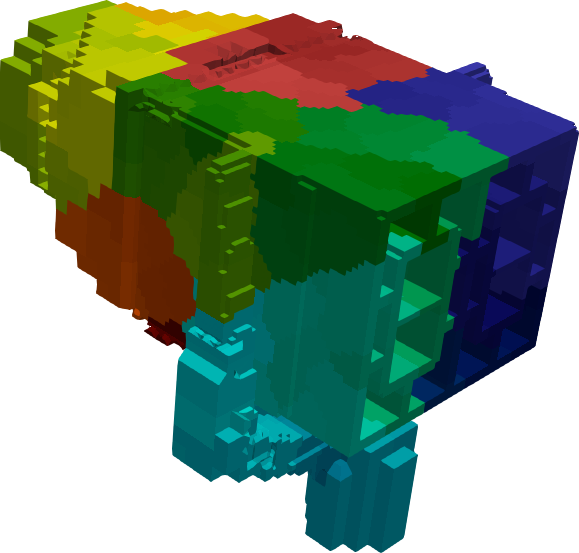} \hspace*{5ex}
  \includegraphics[width=0.45\textwidth]{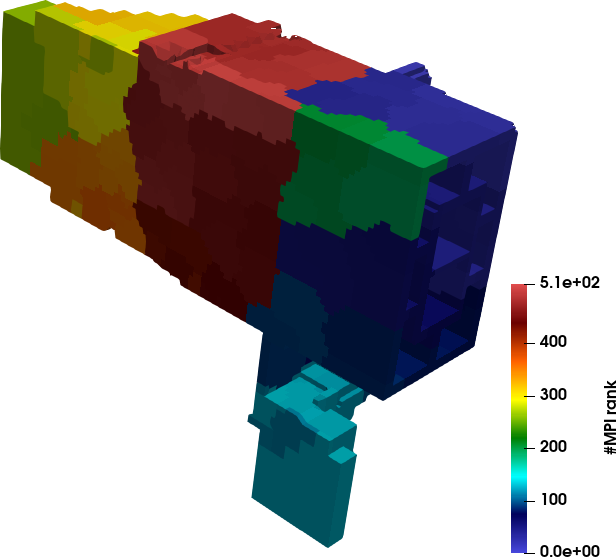}
  \caption{Partition of the extrusion die into 512 subdomains.}\label{pamm-26-koehler-fig:gendie-partition}
\end{figure}

\begin{vchtable}[tpb]
  \centering
  {
    \begin{tabular}{r r | r  r  r | r  r  r | r  r  r}
      \toprule
      \multicolumn{11}{c}{Carreau-type non-Newtonian Stokes on an extrusion die geometry}\\[0.5ex]
       &  & \multicolumn{3}{c|}{2-level FROSch} & \multicolumn{3}{c|}{3-level FROSch} & \multicolumn{3}{c}{4-level FROSch}\\
      \cmidrule(lr){3-5}\cmidrule(lr){6-8}\cmidrule(lr){9-11}
      \#DoFs & \#MPI
      & Time & Effic$_{\mathrm{rel}}$ & \#Avg.
      & Time & Effic$_{\mathrm{rel}}$ & \#Avg.
      & Time & Effic$_{\mathrm{rel}}$ & \#Avg.\\
       & ranks
      & Time & & Krylov.\ it.
      & Time & & Krylov.\ it.
      & Time & & Krylov.\ it.\\
      \midrule
      \multicolumn{11}{c}{\emph{FROSch as preconditioner for FGMRES}}\\[1ex]
      2.10\,M & 2\,048 & 55.9\,s & 65\% & 50.2 &  36.2\,s & 100\% & 66.8 & 39.8\,s  & 91\% & 73.8\\
      13.2\,M  & 4\,000 & 161\,s &  72\% & 55.2 &  81.5\,s & 143\% & 86.3 & 84.6\,s & 138\% & 97.3\\
      \midrule
      \multicolumn{11}{c}{\makecell{\emph{FROSch as coarse solver inside the FEATFLOW geometric multigrid} \\ \emph{(fine problem: $13.2\,\mathrm{M}$ dofs for 2\,048 MPI ranks and $94.1\,\mathrm{M}$ dofs for 4\,000 MPI ranks)}}}\\[2.5ex]
      2.10\,M & 2\,048 & 66.8\,s & 78\% & 11.3 & 51.9\,s & 100\%  & 20.2 & 53.5\,s  & 97\% & 20.6\\
      13.2\,M  & 4\,000 & 176\,s & 95\% & 10.0 & 104\,s & 161\%  & 15.8 & 115\,s  & 145\% & 21.6\\
      \bottomrule

    \end{tabular}}
  \caption{\footnotesize Extrusion die problem: comparison of two-, three-, and four-level monolithic overlapping Schwarz FROSch in two roles within the same outer alternating Picard--Newton iteration. Top: FROSch as preconditioner for FGMRES applied to the linearized system; the \#dof column refers to the problem on which FROSch acts, and total runtimes refer to the full solve. Bottom: FROSch as coarse solver inside the FEATFLOW geometric multigrid solver, applied to the linearized system, for a fine problem of $13.2\,\mathrm{M}$ dofs on $2\,048$ MPI ranks and $94.1\,\mathrm{M}$ dofs on $4\,000$ MPI ranks; the \#DoFs column refers to the coarse problem on which FROSch acts, and total runtimes refer to the full multigrid solve.
    Due to limited computing resources, we could not perform a full scaling study.
    Instead, we report the relative scaled efficiency
$
  {\rm Effic}_{\mathrm{rel}}
  =
  \frac{N/(TP)}{N_0/(T_0P_0)}
  =
  \frac{T_0}{T}\,
  \frac{N}{N_0}\,
  \frac{P_0}{P},
$
where \(T_0\) denotes the fastest baseline runtime obtained with
\(P_0=2048\) MPI ranks for \(N_0=2.1\) million degrees of freedom.  The
quantities \(T\), \(P\), and \(N\) denote the runtime, number of MPI ranks, and
problem size of the run under consideration.
  }
  \label{tab:die-frosch-comparison-2}
\end{vchtable}

The results reported in Table~\ref{tab:die-frosch-comparison-2} are intended as a proof of concept of the monolithic four-level overlapping Schwarz preconditioner on a realistic, geometrically complex configuration. Due to limited computing resources, a full weak or strong scaling study was not feasible. Instead, we compare two configurations -- a $2.10\,\mathrm{M}$-dof problem on $2\,048$ MPI ranks and a $13.2\,\mathrm{M}$-dof problem on $4\,000$ MPI ranks -- and report, in addition to the runtimes, the relative scaled efficiency $\mathrm{Effic}_{\mathrm{rel}}$ defined in the caption, which accounts for the simultaneous change of problem size and rank count. We emphasize that the underlying problem is fully three-dimensional, posed on the complex geometry of an extrusion die, and is solved on an unstructured decomposition of $\Omega$, see Fig.~\ref{pamm-26-koehler-fig:gendie-partition}; the associated load imbalance affects all three preconditioner variants.

Across both roles -- FROSch as preconditioner for FGMRES on the linearized system (top section) and FROSch as coarse solver inside the FEATFLOW geometric multigrid (bottom section) -- the three-level preconditioner outperforms the two-level preconditioner, as expected: once the coarse problem of the two-level method becomes large, replacing its direct factorization by a further overlapping Schwarz step is beneficial. The four-level preconditioner is, in turn, consistently faster than the two-level variant, but slightly slower than the three-level variant. We attribute this to the fact that, for the problem sizes considered here, the coarse problem of the three-level method is not yet large enough to make a further recursion profitable; the additional level then contributes overhead without producing a sufficiently smaller fourth-level coarse problem. This is reflected in the average Krylov iteration counts, which grow mildly with each additional level, in both roles. In the multigrid coarse-solver role (bottom section), the absolute Krylov iteration counts are small -- roughly $10$ for the two-level method and at most about $22$ for the four-level method -- which is a direct consequence of the fact that the multigrid coarse problem is solved only to low accuracy, namely to a relative residual reduction of $10^{-1}$, which is sufficient. This regime, in which a single coarse solve is cheap but is invoked many times across the multigrid cycles and the outer Picard--Newton iterations, is precisely the setting in which multi-level overlapping Schwarz preconditioners are of particular interest, since the per-application overhead of an additional level is amortized only mildly while the underlying coarse factorization is replaced by a scalable approximation.

Several of the reported efficiencies exceed $100\%$. This is a consequence of the definition of the relative scaled efficiency together with the chosen problem sizes. Going from the baseline with $P_{0}=2\,048$ ranks to $P=4\,000$ ranks increases the resources by roughly a factor of two, while the problem grows from $N_{0}=2.10\,\mathrm{M}$ to $N=13.2\,\mathrm{M}$ dofs, i.e.\ by roughly a factor of six.
Since $\mathrm{Effic}_{\mathrm{rel}} = (T_{0}/T)\,(N/N_{0})\,(P_{0}/P)$ and $(N/N_{0})\,(P_{0}/P)\approx 6\cdot\tfrac{1}{2}=3$, a run that takes about three times the baseline runtime $T_{0}$ would already attain an efficiency of $100\%$; any run completing in less than three times $T_{0}$ therefore exceeds $100\%$.
We attribute the fact that the observed efficiencies larger than 100\%
to a load imbalance in the problem with 2.10M unknowns.
This is indicated by the size of the largest subdomain: For 2.10M unknowns and $2048$ ranks, the largest first-level subdomain has approximately $11\,\mathrm{K}$ dofs, whereas for 13.2M unknowns and $4000$ ranks the largest subdomain has approximately $21\,\mathrm{K}$ unknowns.

A weak scaling study on larger supercomputers, in which the second-level coarse problem of the three-level method becomes the bottleneck, is required to reach the regime in which the four-level method is superior.

\begin{acknowledgement}
  This work is funded by the German Federal Ministry of Research, Technology and Space (BMFTR, formerly BMBF) under grant no.~16ME0708 as part of the SCALEXA initiative on exascale computing.
  Funded by the European Union -- NextGenerationEU.

  The authors acknowledge computing time on the compute cluster of the Faculty of Mathematics and Computer Science of Technische Universität Bergakademie Freiberg, operated by the computing center (URZ) and funded by the Deutsche Forschungsgemeinschaft (DFG) under DFG grant number 397252409 (\url{https://gepris.dfg.de/gepris/projekt/397252409}).

  We would like to thank Peter Zajac and Maximilian Esser from TU Dortmund for their part in the implementation of the FEATFLOW/FROSch interface and in the numerical examples, and Otto Mierka for the meshes.

  During the preparation of this work the generative AI tools Claude were used to correct grammar, spelling, and to improve style of writing. After using these tools, all authors reviewed and edited the content as needed and take full responsibility for the content of the
  publication.
\end{acknowledgement}

\vspace{\baselineskip}
\bibliographystyle{pamm-doi}
\bibliography{all_refs}

\end{document}